\documentclass{article}
\usepackage{graphicx}
\usepackage{amssymb, amsthm, amsmath, amsfonts}
\usepackage{palatino}
\usepackage[T1]{fontenc}
\usepackage[sc]{mathpazo}
\usepackage{inconsolata}
\usepackage{titlesec}
\usepackage{xcolor}
\usepackage[sort,numbers]{natbib}
\usepackage{mathdots}
\usepackage{enumerate}
\usepackage{caption}

\usepackage{lineno}
\newcommand*\patchAmsMathEnvironmentForLineno[1]{%
  \expandafter\let\csname old#1\expandafter\endcsname\csname #1\endcsname
  \expandafter\let\csname oldend#1\expandafter\endcsname\csname end#1\endcsname
  \renewenvironment{#1}%
     {\linenomath\csname old#1\endcsname}%
     {\csname oldend#1\endcsname\endlinenomath}}% 
\newcommand*\patchBothAmsMathEnvironmentsForLineno[1]{%https://www.overleaf.com/19542727zwszmncxjmbg
  \patchAmsMathEnvironmentForLineno{#1}%
  \patchAmsMathEnvironmentForLineno{#1*}}%
\AtBeginDocument{%
\patchBothAmsMathEnvironmentsForLineno{equation}%
\patchBothAmsMathEnvironmentsForLineno{align}%
\patchBothAmsMathEnvironmentsForLineno{flalign}%
\patchBothAmsMathEnvironmentsForLineno{alignat}%
\patchBothAmsMathEnvironmentsForLineno{gather}%
\patchBothAmsMathEnvironmentsForLineno{multline}%
}
\definecolor{todocolor}{RGB}{205,235,139}
\definecolor{todo-idea}{RGB}{120,180,255}
\definecolor{todo-error}{RGB}{208,31,60}
\definecolor{todo-question}{RGB}{255,255,136}

\usepackage[colorinlistoftodos, color=todocolor, textsize=small]{todonotes}

\usepackage{enumerate}

\usepackage[margin=1.25in]{geometry}

\theoremstyle{plain}
\newtheorem{theorem}{Theorem}%[section]
\newtheorem{corollary}[theorem]{Corollary}

\newtheorem{porism}[theorem]{Porism}

\newcommand{\ds}{\displaystyle}
\newcommand{\lc}[2]{{#1}^{c_{#2}}}

\newcommand{\eu}[2]{\genfrac<>{0pt}{}{#1}{#2}}

\title{
    Colored Multipermutations and a Combinatorial Generalization of Worpitzky's Identity
}

\author{John Engbers\thanks{john.engbers@marquette.edu; Department of Mathematical and Statistical Sciences, Marquette University, Milwaukee, WI 53201. Research supported by the Simons Foundation grant 524418.}  
\and Jay Pantone\thanks{jay.pantone@marquette.edu; Department of Mathematical and Statistical Sciences, Marquette University, Milwaukee, WI 53201.}  
\and Christopher Stocker\thanks{christopher.stocker@marquette.edu; Department of Mathematical and Statistical Sciences, Marquette University, Milwaukee, WI 53201.}}

\date{}

\begin{document}

\maketitle

\begin{abstract}
    \setlength{\parskip}{9pt plus 3pt minus 1pt}
    \setlength{\parindent}{0pt}
    
    Worpitzky's identity~\cite{worpitzky:identity} expresses $n^p$ in terms of the Eulerian numbers and binomial coefficients:
    \[
        n^p = \sum_{i=0}^{p-1}\eu{p}{i}\binom{n+i}{p}.
    \]
    Pita-Ruiz~\cite{pita-ruiz:generalization-eulerian,pita-ruiz:generalized-eulerian-and-apps} recently defined numbers $A_{a,b,r}(p,i)$ implicitly to satisfy a generalized Worpitzky identity
    \[
        \binom{an+b}{r}^p = \sum_{i=0}^{rp} A_{a,b,r}(p,i) \binom{n+rp-i}{rp},
    \]
    and asked whether there is a combinatorial interpretation of the numbers $A_{a,b,r}(p,i)$.
    
    We provide such a combinatorial interpretation by defining a notion of descents in colored multipermutations, and then proving that $A_{a,b,r}(p,i)$ is equal to the number of colored multipermutations of $\{1^r, 2^r, \ldots, p^r\}$ with $a$ colors and $i$ weak descents. We use this to give combinatorial proofs of several identities involving $A_{a,b,r}(p,i)$, including the aforementioned generalized Worpitzky identity.
\end{abstract}

%====================================================================================

\section{Generalized Eulerian Numbers}

The famous \emph{Eulerian numbers} $\ds\eu{p}{i}$ appear first in Euler's 1755 manuscript \emph{Institutiones calculi differentialis}~\cite{euler:foundations} and are now known to count a variety of combinatorial objects. The Eulerian numbers can be defined by the recurrence
\[
    \eu{p}{i} = (p-i)\eu{p-1}{i-1} + (i+1)\eu{p-1}{i}
\]
for $1 \leq i \leq p-1$, with $\ds\eu{p}{0} = \ds\eu{p}{p-1} = 1$ for $p \geq 1$.

Among the many identities involving the Eulerian numbers is Worpitzky's identity~\cite{worpitzky:identity},
\[
    n^p = \sum_{i=0}^{p-1} \eu{p}{i}\binom{n+i}{p},
\]
which demonstrates that the Eulerian numbers are the necessary coefficients to write the powers $n^p$ using a basis of binomial coefficients. Moreover, an explicit formula is known:
\begin{equation}\label{eq-Eul}
    \eu{p}{i} = \sum_{j=0}^{i+1} (-1)^j \binom{p+1}{j} (i-j+1)^p.
\end{equation}
We refer interested readers to Petersen's excellent book \emph{Eulerian Numbers}~\cite{petersen:eulerian-numbers}.

In a recent pair of papers, Pita-Ruiz~\cite{pita-ruiz:generalization-eulerian,pita-ruiz:generalized-eulerian-and-apps} defines a family of generalized Eulerian numbers $A_{a,b,r}(p,i)$ and proves an incredible number of identities, explicit formulas, and recurrences. While there are many different generalizations of the Eulerian numbers (see~\cite{rzadkowski:generalization-euler} and the references therein), it is notable for the purposes of this paper that Pita-Ruiz's definition is not combinatorially-based, but rather an implicit definition  to satisfy a generalization of Worpitzky's identity:
\begin{equation}
\label{equation:def-of-A}
    \binom{an+b}{r}^p = \sum_{i=0}^{rp} A_{a,b,r}(p,i) \binom{n+rp-i}{rp}.
\end{equation}
Note that when $a=1$ and $b=0$, the left side involves powers of binomials $\ds\binom{n}{r}$, and analytical proofs of equations involving summations of these powers, and some generalizations, were given by Alzer and Prodinger~\cite{alzer:combinatorial-identities} with corresponding combinatorial proofs given by Engbers and Stocker~\cite{engbers:two-comb-proofs, engbers:comb-proofs-Alzer-Prodinger}.

Pita-Ruiz~\cite{pita-ruiz:generalization-eulerian,pita-ruiz:generalized-eulerian-and-apps} derives a summation formula
\[
    A_{a,b,r}(p,i) = \sum_{j=0}^{i} (-1)^j \binom{rp+1}{j} \binom{a(i-j)+b}{r}^p
\]
that bears some resemblance to the summation formula for the ordinary Eulerian numbers in \eqref{eq-Eul}, and asks whether there is a combinatorial interpretation of these numbers $A_{a,b,r}(p,i)$.

One of the many interpretations of the ordinary Eulerian numbers is that $\ds\eu{p}{i}$ is equal to the number of permutations of $\{1, \ldots, p\}$ with $i$ descents\footnote{A descent in a permutation $\pi = \pi(1)\cdots\pi(p)$ is a pair of consecutive entries $\pi(i)\pi(i+1)$ with $\pi(i) > \pi(i+1)$.}. In the following sections, we provide a combinatorial interpretation of $A_{a,b,r}(p,i)$.  In particular, we define a type of colored multipermutation with a generalized notion of a descent such that:

\emph{$A_{a,b,r}(p,i)$ is the number of colored multipermutations on $p$ symbols with multiplicity $r$, with $a$ colors, and with $i$ weak descents.}

We remark that $b$ is a parameter that varies the definition of a weak descent.

In Section \ref{sec:coloredmulti}, we describe the colored multipermutations that we use in our combinatorial interpretation, and then utilize these in Section \ref{sec:b=0} to provide the interpretation when $b=0$.  The extension to the case $0<b<a$, followed by remarks about $b \geq a$, appears in Section \ref{sec:b>0}.  Finally, in Section \ref{sec:Worpitzky} we use the combinatorial interpretation to provide a combinatorial proof of the generalized Worpitzky identity given in  Equation~\eqref{equation:def-of-A}.

\section{Colored Multipermutations}\label{sec:coloredmulti}

For positive integers $p$ and $r$, let $\{1^r, 2^r, \ldots, p^r\}$ denote the multiset that contains $r$ copies of each element of $\{1, 2, \ldots, p\}$. A linear ordering of the elements of $\{1^r, 2^r, \ldots, p^r\}$ is called a \emph{multipermutation}; for example, $311232$ is a multipermutation of $\{1^2, 2^2, 3^2\}$.

A \emph{colored multipermutation with $a$ colors} is a multipermutation in which entry has been assigned one of $a$ colors. There are no restrictions to how these colors may be assigned. We represent the colors $c_1, \ldots, c_a$ as superscripts, e.g.,
\[
    \lc{3}{1}\lc{1}{5}\lc{1}{2}\lc{2}{7}\lc{3}{2}\lc{2}{2}
\]
can be considered a colored multipermutation of $\{1^2, 2^2,3^2\}$ with $9$ colors, even though we have only used four of the available colors. From these definitions, it can be readily seen that the number of colored multipermutations of $\{1^r, 2^r, \ldots, p^r\}$ with $a$ colors is
\[
    a^{rp} \frac{(rp)!}{(r!)^p}.
\]

In order to define the notion of a descent in a colored multipermutation, we define a natural ordering on the set of colored entries, which is the lexicographic order first on color, then on value:
\[
    \lc{1}{1} < \lc{2}{1} < \cdots < \lc{p}{1} < \lc{1}{2} < \lc{2}{2} < \cdots < \lc{p}{2} < \quad\cdots\quad < \lc{1}{a} < \lc{2}{a} < \cdots < \lc{p}{a}.
\]
A \emph{non-terminal descent} in a colored multipermutation is a pair of adjacent entries $\lc{i_1}{j_1}\lc{i_2}{j_2}$ such that $\lc{i_1}{j_1} > \lc{i_2}{j_2}$ according to the order defined above. Further, a \emph{non-terminal weak descent} (which could alternatively be called a \emph{non-terminal non-ascent}) is a pair of adjacent entries $\lc{i_1}{j_1}\lc{i_2}{j_2}$ such that $\lc{i_1}{j_1} \geq \lc{i_2}{j_2}$.

We now introduce one final parameter $b$ that determines whether or not there is a \emph{terminal descent} after the final entry of a colored multipermutation. Let $0 \leq b < a$ be an integer. The final entry $\lc{i}{j}$ is considered to be itself a terminal descent if its color has value greater than $b$, i.e., $j > b$. When $b=0$, for example, the final entry is always a terminal descent, while when $b = a-1$, the final entry is a terminal descent only if its color is $c_a$. The set of \emph{weak descents} consists of the set of non-terminal weak descents, together with the terminal descent if there is one. We call a colored multipermutation of $\{1^r, \ldots, p^r\}$ with $a$ colors ``$(p,r,a,b)$-type'' when its descents should be counted with respect to the parameter $b$.

The ideas of adding colors or multiplicity to permutations are not new---they have been studied in myriad combinations, including with regard to their descents (which have also been defined in many ways); see, e.g.,~\cite{park:r-multipermutations,moynihan:colored-eulerian-descent-algebra,chen:labeled-partitions-colored-perms, bagno:counting-descent-pairs}.

In the following sections, we prove that the family of numbers $A_{a,b,r}(p,i)$ defined implicitly by Pita-Ruiz~\cite{pita-ruiz:generalization-eulerian,pita-ruiz:generalized-eulerian-and-apps} to satisfy Equation~(\ref{equation:def-of-A}) has in fact a very nice combinatorial interpretation: {\em $A_{a,b,r}(p,i)$ is the number of $(p,r,a,b)$-type colored multipermutations with $i$ weak descents}. This  interpretation allows us to easily provide combinatorial proofs for several identities of Pita-Ruiz~\cite{pita-ruiz:generalization-eulerian,pita-ruiz:generalized-eulerian-and-apps} involving $A_{a,b,r}(p,i)$.

\section{The $b=0$ Case}\label{sec:b=0}

We start with a thorough explanation of the case when $b=0$, from which the $0 < b < a$ case follows with only minor adaptations. When $b=0$, the last entry is always considered to be a weak descent. Recall that we call this the terminal weak descent and call all other weak descents non-terminal. The following theorem provides a formula for the number of colored multipermutations with a fixed number of non-terminal weak descents, and the simple corollary following it recovers combinatorially the formula of Pita-Ruiz~\cite{pita-ruiz:generalization-eulerian,pita-ruiz:generalized-eulerian-and-apps}.

\begin{theorem}
    \label{theorem:main-b-equals-0-case}
    The number of colored multipermutations of $\{1^r, \ldots, p^r\}$ with $a$ colors and $i$ non-terminal weak descents is
    \[
        \sum_{j=0}^{i+1 } (-1)^j \binom{rp+1}{j} \binom{a(i-j+1)}{r}^p.
    \]
\end{theorem}
\begin{proof}
    
    Fix positive integers $r$, $p$, $a$, and $i$. We first define a \emph{segmented colored multipermutation} (SCM) with $i$ segments to be a colored multipermutation split into $i$ contiguous parts by $i-1$ dividing lines with the property that each segment contains no non-terminal weak descents. We call the dividing lines \emph{walls}. It is permissible for some segments to be empty, so that two or more walls are adjacent. For example, 
    \[
        \lc{2}{1}\lc{4}{1}\lc{1}{3}\lc{3}{3}
        \,\Bigg|\,
        \lc{1}{1}\lc{2}{1}\lc{4}{1}
        \,\Bigg|\,
        \lc{2}{2}\lc{5}{2}\lc{5}{3}
        \,\Bigg|\,\,\Bigg|\,
        \lc{1}{1}\lc{5}{1}\lc{4}{2}\lc{3}{3}
        \,\Bigg|\,
        \lc{3}{3}
    \]
    is an SCM with six segments (and thus five walls). %Observe that each of the six segments contains no non-terminal weak descents when taken in isolation and that the fourth segment is empty. Note also 
    Note that the underlying colored multipermutation contains only three non-terminal weak descents, found between the pairs of entries on either side of the first, third, and fifth walls. 
    
    We call a wall in an SCM \emph{extraneous} if either of two conditions hold:
    \begin{enumerate}[(a)]
        \item the segments to the left and right are nonempty, and the deletion of the wall creates a larger segment with no non-terminal weak descents, or
        \item the segment to the left is empty.
    \end{enumerate}
    
    The first wall in the example SCM given above is not extraneous because its removal creates the segment
    \[
        \lc{2}{1}\lc{4}{1}\lc{1}{3}\lc{3}{3}
        \lc{1}{1}\lc{2}{1}\lc{4}{1}
    \]
    which contains the weak descent $\lc{3}{3}\lc{1}{1}$. On the other hand, the second wall is extraneous because its removal creates the segment
    \[
        \lc{1}{1}\lc{2}{1}\lc{4}{1}
        \lc{2}{2}\lc{5}{2}\lc{5}{3}
    \]
    which contains no weak descents. The fourth wall is extraneous because the segment to its left is empty, and the third and fifth walls are not extraneous.
    
    There is an obvious bijection %$\alpha$ 
    between the set of $(p,r,a,0)$-type colored multipermutations with $i$ non-terminal weak descents and the set of SCMs with $i+1$ segments (thus $i$ walls) and no extraneous walls and whose underlying colored multipermutation has type $(p,r,a,0)$. This bijection is easily described as inserting or removing walls between all non-terminal weak descents. As a consequence, we focus now on the enumeration of SCMs with $i+1$ segments and no extraneous walls. With the aim of eventually applying the inclusion-exclusion principle, we start by calculating the enumeration of SCMs with $i+1$ segments and no restriction on the number of extraneous walls.
    
    An SCM with $i+1$ segments whose underlying colored multipermutation has type $(p,r,a,0)$ can be constructed bijectively with the following procedure. Imagine that there are $i+1$ columns, each with $a$ bins, arranged from left-to-right as in Figure~\ref{figure:bins-1}.
    
    \begin{figure}[h]
        \begin{center}
            \begin{tikzpicture}
                \draw[thick,darkgray,rounded corners=0.5ex] (0,0) rectangle (1,0.5);
                \draw[thick,darkgray,rounded corners=0.5ex] (0,0.75) rectangle (1,1.25);
                \node at (0.5, 1.85) {$\vdots$};
                \draw[thick,darkgray,rounded corners=0.5ex] (0,2.25) rectangle (1,2.75);
                
                \begin{scope}[shift={(1.5,0)}]
                    \draw[thick,darkgray,rounded corners=0.5ex] (0,0) rectangle (1,0.5);
                    \draw[thick,darkgray,rounded corners=0.5ex] (0,0.75) rectangle (1,1.25);
                    \node at (0.5, 1.85) {$\vdots$};
                    \draw[thick,darkgray,rounded corners=0.5ex] (0,2.25) rectangle (1,2.75);
                \end{scope}
                
                \begin{scope}[shift={(4.5,0)}]
                    \draw[thick,darkgray,rounded corners=0.5ex] (0,0) rectangle (1,0.5);
                    \draw[thick,darkgray,rounded corners=0.5ex] (0,0.75) rectangle (1,1.25);
                    \node at (0.5, 1.85) {$\vdots$};
                    \draw[thick,darkgray,rounded corners=0.5ex] (0,2.25) rectangle (1,2.75);
                \end{scope}
                
                \node at (3.5, 0.2) {$\cdots$};
                \node at (3.5, 0.95) {$\cdots$};
                % \node at (3.5, 1.85) {$\iddots$};
                \node at (3.5, 2.45) {$\cdots$};
                
                \node[darkgray] at (0.5, -0.25) {\scriptsize \texttt{column $1$}};
                \node[darkgray] at (2, -0.25) {\scriptsize \texttt{column $2$}};
                \node[darkgray] at (5, -0.25) {\scriptsize \texttt{column $i+1$}};
                
                % \node[darkgray, rotate=270] at (-0.25, 0.5) {\scriptsize \texttt{color $1$}};
                % \node[darkgray] at (-0.25, 0.2) {\scriptsize \texttt{$c_1$}};
                \node[darkgray] at (-0.6, 0.25) {\scriptsize \texttt{color $1$}};
                \node[darkgray] at (-0.6, 1) {\scriptsize \texttt{color $2$}};
                \node[darkgray] at (-0.6, 2.5) {\scriptsize \texttt{color $a$}};
                
                \draw[gray] (1.25, -0.25) -- (1.25, 2.9);
                \draw[gray] (2.75, -0.25) -- (2.75, 2.9);
                \draw[gray] (4.25, -0.25) -- (4.25, 2.9);
            \end{tikzpicture}
        \end{center}
        \vspace*{-0.6cm}
        \captionsetup{width=.8\linewidth}
        \caption{An arrangement of $a(i+1)$ bins into $i+1$ columns and $a$ rows.}
        \label{figure:bins-1}
    \end{figure}
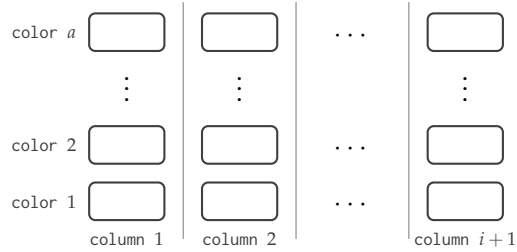
    
    Then, distribute $r$ copies of each of the symbols $1, 2, \ldots, p$ into these $a(i+1)$ bins such that each bin contains at most one copy of each symbol (but different symbols are allowed in the same bin). This can be done in $\binom{a(i+1)}{r}^p$ ways. An SCM is then formed from this choice by arranging the entries in each column in the unique manner with no non-terminal weak descents and concatenating the resulting words from left-to-right while placing a wall between the entries from each column. Figure~\ref{figure:bins-assignment} shows the assignment of entries to bins for $r = 3$ and $p = 5$ that results in the SCM used as an example above:
    \[
        \lc{2}{1}\lc{4}{1}\lc{1}{3}\lc{3}{3}
        \,\Bigg|\,
        \lc{1}{1}\lc{2}{1}\lc{4}{1}
        \,\Bigg|\,
        \lc{2}{2}\lc{5}{2}\lc{5}{3}
        \,\Bigg|\,\,\Bigg|\,
        \lc{1}{1}\lc{5}{1}\lc{4}{2}\lc{3}{3}
        \,\Bigg|\,
        \lc{3}{3}.
    \]
    This shows that the number of SCMs with $i+1$ segments and no constraints on the number of extraneous walls is $\binom{a(i+1)}{r}^p$.
    
    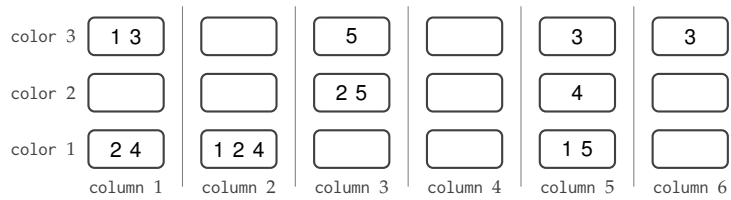
\begin{figure}[h]
        \begin{center}
            \begin{tikzpicture}
                \draw[thick,darkgray,rounded corners=0.5ex] (0,0) rectangle (1,0.5);
                \draw[thick,darkgray,rounded corners=0.5ex] (0,0.75) rectangle (1,1.25);
                \draw[thick,darkgray,rounded corners=0.5ex] (0,1.5) rectangle (1,2);
                \node[darkgray] at (0.5, -0.25) {\scriptsize \texttt{column $1$}};
                \draw[gray] (1.25, -0.25) -- (1.25, 2.15);
                \node at (0.5, 0.25) {\footnotesize \textsf{2}\;\;\!\textsf{4}};
                \node at (0.5, 1) {\footnotesize };
                \node at (0.5, 1.75) {\footnotesize \textsf{1}\;\;\!\textsf{3}};
                
                \begin{scope}[shift={(1.5,0)}]
                    \draw[thick,darkgray,rounded corners=0.5ex] (0,0) rectangle (1,0.5);
                    \draw[thick,darkgray,rounded corners=0.5ex] (0,0.75) rectangle (1,1.25);
                    \draw[thick,darkgray,rounded corners=0.5ex] (0,1.5) rectangle (1,2);
                    \node[darkgray] at (0.5, -0.25) {\scriptsize \texttt{column $2$}};
                    \draw[gray] (1.25, -0.25) -- (1.25, 2.15);
                    \node at (0.5, 0.25) {\footnotesize \textsf{1}\;\;\!\textsf{2}\;\;\!\textsf{4}};
                \end{scope}
                
                \begin{scope}[shift={(3.0,0)}]
                    \draw[thick,darkgray,rounded corners=0.5ex] (0,0) rectangle (1,0.5);
                    \draw[thick,darkgray,rounded corners=0.5ex] (0,0.75) rectangle (1,1.25);
                    \draw[thick,darkgray,rounded corners=0.5ex] (0,1.5) rectangle (1,2);
                    \node[darkgray] at (0.5, -0.25) {\scriptsize \texttt{column $3$}};
                    \draw[gray] (1.25, -0.25) -- (1.25, 2.15);
                    \node at (0.5, 1) {\footnotesize \textsf{2}\;\;\!\textsf{5}};
                    \node at (0.5, 1.75) {\footnotesize \textsf{5}};
                \end{scope}
                
                \begin{scope}[shift={(4.5,0)}]
                    \draw[thick,darkgray,rounded corners=0.5ex] (0,0) rectangle (1,0.5);
                    \draw[thick,darkgray,rounded corners=0.5ex] (0,0.75) rectangle (1,1.25);
                    \draw[thick,darkgray,rounded corners=0.5ex] (0,1.5) rectangle (1,2);
                    \node[darkgray] at (0.5, -0.25) {\scriptsize \texttt{column $4$}};
                    \draw[gray] (1.25, -0.25) -- (1.25, 2.15);
                \end{scope}
                
                \begin{scope}[shift={(6.0,0)}]
                    \draw[thick,darkgray,rounded corners=0.5ex] (0,0) rectangle (1,0.5);
                    \draw[thick,darkgray,rounded corners=0.5ex] (0,0.75) rectangle (1,1.25);
                    \draw[thick,darkgray,rounded corners=0.5ex] (0,1.5) rectangle (1,2);
                    \node[darkgray] at (0.5, -0.25) {\scriptsize \texttt{column $5$}};
                    \draw[gray] (1.25, -0.25) -- (1.25, 2.15);
                    \node at (0.5, 0.25) {\footnotesize \textsf{1}\;\;\!\textsf{5}};
                    \node at (0.5, 1) {\footnotesize \textsf{4}};
                    \node at (0.5, 1.75) {\footnotesize \textsf{3}};
                \end{scope}
                
                \begin{scope}[shift={(7.5,0)}]
                    \draw[thick,darkgray,rounded corners=0.5ex] (0,0) rectangle (1,0.5);
                    \draw[thick,darkgray,rounded corners=0.5ex] (0,0.75) rectangle (1,1.25);
                    \draw[thick,darkgray,rounded corners=0.5ex] (0,1.5) rectangle (1,2);
                    \node[darkgray] at (0.5, -0.25) {\scriptsize \texttt{column $6$}};
                    \node at (0.5, 1.75) {\footnotesize \textsf{3}};
                \end{scope}
                
                \node[darkgray] at (-0.6, 0.25) {\scriptsize \texttt{color $1$}};
                \node[darkgray] at (-0.6, 1) {\scriptsize \texttt{color $2$}};
                \node[darkgray] at (-0.6, 1.75) {\scriptsize \texttt{color $3$}};
                
            \end{tikzpicture}
        \end{center}
        \vspace*{-0.6cm}
        \captionsetup{width=.8\linewidth}
        \caption{An assignment of the elements of $\{1^3,2^3,3^3,4^3,5^3\}$ into an arrangement of bins with $6$ columns and $3$ rows.}
        \label{figure:bins-assignment}
    \end{figure}
    
    Define the \emph{position} of a wall in an SCM to be the number of entries to its left. For example, in an SCM whose underlying colored multipermutation has $rp$ entries, a wall to the left of all entries is in position $0$ and a wall to the right of all entries is in position $rp$. Two walls can be in the same position if the segment between them is empty.
    
    For any set $S \subseteq \{0, 1, \ldots, rp\}$ let $A_S$ denote the set of SCMs of $(p, r, a, 0)$ type with $i+1$ segments such that there is at least one extraneous wall in each position in $S$. There are permitted to be extraneous walls in other positions as well. The set $A_{\emptyset}$ consists of all SCMs of $(p, r, a, 0)$ type with $i+1$ segments. Now set
    \begin{align*}
        Q &= A_{\emptyset} \smallsetminus \bigcup_{\substack{S \subseteq \{0, 1, \ldots, rp\}\\S \neq \emptyset}} A_S\\
        &= \bigcap_{\substack{S \subseteq \{0, 1, \ldots, rp\}\\S \neq \emptyset}} \overline{A_S}.
    \end{align*}
    The set $Q$ consists of the SCMs of $(p, r, a, 0)$ type with $i+1$ segments and no extraneous walls. We observed earlier that $Q$ is in bijection with the set of objects we wish to count, the colored multipermutations with exactly $i$ non-terminal weak descents. By the inclusion-exclusion principle, the size of $Q$ can be calculated by
    \[
        |Q| = \sum_{S \subseteq \{0, 1, \ldots rp\}} (-1)^{|S|}|A_S|,
    \]
    and so it remains only to determine $|A_S|$. For this we define one final bijection, inspired by the proof of Theorem 1.11 in B\'ona's \emph{Combinatorics of Permutations}~\cite{bona:combinatorics-of-permutations}.
    
    For any set $S \subseteq \{0, 1, \ldots, rp\}$, define a function $g_S$ whose domain is $A_S$ and whose codomain is the set of $(p, r, a, 0)$-type SCMs with $i+1-|S|$ segments as follows. For $\alpha$ in the domain of $g_S$, form $\beta$ by deleting one wall in each of the positions in $S$. Since $\alpha$ must have at least one extraneous wall in each of the positions in $S$ (by the definition of $A_S$), the result $\beta$ is a valid SCM with $i-|S|$ walls and thus $i+1-|S|$ segments. The inverse map takes as input an SCM $\beta$ with $i-|S|$ walls and adds one wall in each position of $S$ (if a wall is already present, the added wall is placed on the right of the existing wall). Each added wall is then guaranteed to be extraneous, and the resulting SCM has $i+1$ segments.
    
    We have already calculated earlier in this proof that the number of SCMs with $\ell$ segments is $\binom{a\ell}{r}^p$, and thus the bijection above demonstrates that
    \[
        |A_S| = \ds\binom{a(i+1-|S|)}{r}^p.
    \]
    
    Plugging this into the result of the inclusion-exclusion above yields
    \begin{align*}
        |Q| &= \sum_{S \subseteq \{0, 1, \ldots rp\}} (-1)^{|S|}|A_S|\\
        % &= \sum_{S \subseteq \{0, 1, \ldots rp\}} (-1)^{|S|}\binom{a(i+1-|S|)}{r}^p\\
 %       &= \sum_{j=0}^{rp+1} (-1)^j\binom{rp+1}{j}\binom{a(i-j+1)}{r}^p\\
        &= \sum_{j=0}^{i+1} (-1)^j\binom{rp+1}{j}\binom{a(i-j+1)}{r}^p,
    \end{align*}
    which completes the proof.
\end{proof}

\begin{corollary}
    The number of colored multipermutations of $\{1^r, \ldots, p^r\}$ with $a$ colors and $i$ weak descents is 
    \[
        \sum_{j=0}^{i} (-1)^j \binom{rp+1}{j} \binom{a(i-j)}{r}^p,
    \]
    and this quantity is equal to $A_{a,0,r}(p,i)$ as defined by Pita-Ruiz~\cite{pita-ruiz:generalization-eulerian,pita-ruiz:generalized-eulerian-and-apps}.
\end{corollary}
\begin{proof}
    The statement of the corollary differs from the statement of Theorem~\ref{theorem:main-b-equals-0-case} in three ways: ``$i$ non-terminal weak descents'' has been replaced by ``$i$ weak descents'', the upper bound of summation has changed from $i+1$ to $i$, and the upper term in the second binomial coefficients has changed from $a(i-j+1)$ to $a(i-j)$.
    
    The quantity given in Theorem~\ref{theorem:main-b-equals-0-case} counts the number of $(p, r, a, 0)$ colored multipermutations with $i$ non-terminal weak descents, and since the final entry of such a colored multipermutation is always a descent, this can be accounted for by shifting the $i$-index by one, giving
    \[
        \sum_{j=0}^{i} (-1)^j \binom{rp+1}{j} \binom{a(i-j)}{r}^p = A_{a,0,r}(p,i)
    \]
    as claimed.
\end{proof}

\section{The $0 < b < a$ Case}\label{sec:b>0}

We now adapt the proof from the previous section to allow for the case where $b$ is constrained by $0 < b < a$.

\begin{theorem}
    \label{theorem:main-b-neq-0}
    Suppose $0 < b < a$. The number of colored multipermutations of $\{1^r, \ldots, p^r\}$ with $a$ colors and $i$ weak descents, with the consideration that the final position is considered a weak descent if and only if the color of the last entry has value strictly greater than $b$, is
    \[
       \sum_{j=0}^{i} (-1)^j \binom{rp+1}{j} \binom{a(i-j)+b}{r}^p,
    \]
    and this quantity is equal to $A_{a,b,r}(p,i)$ as defined by Pita-Ruiz~\cite{pita-ruiz:generalization-eulerian,pita-ruiz:generalized-eulerian-and-apps}.
\end{theorem}
\begin{proof}
    As in the proof of Theorem~\ref{theorem:main-b-equals-0-case}, we start by defining an auxiliary object. A $b$-segmented colored multipermutation ($b$-SCM) with $i$ segments is as before, a colored multipermutation split into $i$ contiguous parts by $i-1$ dividing lines with the property that each segment contains no weak descents, together with the additional constraint that the color of all entries in the final segment must have color value at most $b$. The example
     \[
        \lc{2}{1}\lc{4}{1}\lc{1}{3}\lc{3}{3}
        \,\Bigg|\,
        \lc{1}{1}\lc{2}{1}\lc{4}{1}
        \,\Bigg|\,
        \lc{2}{2}\lc{5}{2}\lc{5}{3}
        \,\Bigg|\,\,\Bigg|\,
        \lc{1}{1}\lc{5}{1}\lc{4}{2}\lc{3}{3}
        \,\Bigg|\,
        \lc{3}{3}
    \]
    used previously, for which the color of the final entry has value $3$, is a $3$-SCM, but is neither a $1$-SCM nor a $2$-SCM. A wall is called \emph{extraneous} if one of three conditions hold (the first two being identical to those in the proof of Theorem~\ref{theorem:main-b-equals-0-case}):
    \begin{enumerate}[(a)]
        \item the segments to the left and right are nonempty, and the deletion of the wall creates a larger segment with no weak descents,
        \item the segment to the left is empty, or
        \item the wall is the final wall, the segment to its right is empty, and the segment to its left contains only entries whose color value is at most $b$. 
    \end{enumerate}
    If the example above contained a wall in the final position, that wall would be extraneous when $b \geq 3$ (because in this case there is \emph{not} a descent in the final position) and non-extraneous when $b=1,2$ (because in this case there is).
    
    We claim at this point that the set of $(p,r,a,b)$-colored multipermutations with $i$ weak descents is equinumerous to the set of $b$-SCMs with $i+1$ segments, no extraneous walls, and whose underlying colored multipermutation has $(p,r,a,b)$ type. The bijection mimics that from the proof of Theorem~\ref{theorem:main-b-equals-0-case}---a colored multipermutation is mapped to a $b$-SCMs by adding walls between each weak descent and adding a wall after the last entry if it has color value strictly greater than $b$; a $b$-SCM with $i+1$ segments and no extraneous walls is mapped to a colored multipermutation by simply removing the walls.
    
    To complete the proof, we must enumerate the set $Q$ of $b$-SCMs with $i+1$ segments and no extraneous walls. Turning again to the inclusion-exclusion principle, for any set $S \subseteq \{0,1,\ldots,rp\}$ define $A_S$ to be the set of $b$-SCMs of $(p,r,a,b)$ type with $i+1$ segments such that there is at least one extraneous wall in each position in $S$. Then,
    \[
        |Q| = \sum_{S \subseteq \{0,1,\ldots,rp\}} (-1)^{|S|} |A_S|.
    \]
    We can determine $|A_S|$ as before by observing that $A_S$ is in bijection with the set of $b$-SCMs of $(p,r,a,b)$ type with $i+1-|S|$ segments, the bijection again involving just adding/removing the extraneous walls.
    
    To construct a $b$-SCM of $(p,r,a,b)$ type with $i+1-|S|$ segments, start by considering $i+1-|S|$ columns of bins, the first $i-|S|$ each containing $a$ bins and the final containing only $b$ bins, as shown in Figure~\ref{figure:bins-assignment-b-neq-0}. Out of the $a(i-|S|)+b$ bins, we must choose $r$ in which to place $1$, $r$ in which to place $2$, etc., yielding a total of
    \[
        \binom{a(i-|S|)+b}{r}^p
    \]
    placements. Each placement corresponds uniquely to a $b$-SCM with $i+1-|S|$ segments by assigning colors to each entry according to which bin they are in, ordering the entries in each column in the unique way with no weak descents, and lining them up from left to right with a wall between each grouping. 
    
    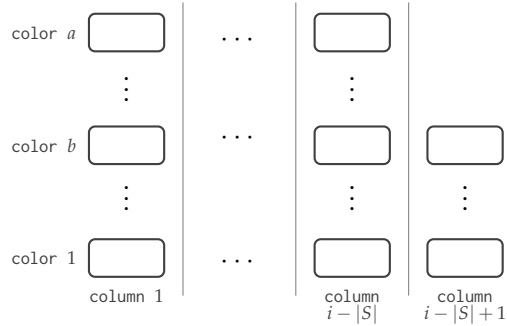
\begin{figure}[ht]
        \begin{center}
            \begin{tikzpicture}
                \draw[thick,darkgray,rounded corners=0.5ex] (0,0) rectangle (1,0.5);
                \node at (0.5, 1.15) {$\vdots$};
                \draw[thick,darkgray,rounded corners=0.5ex] (0,1.5) rectangle (1,2);
                \node at (0.5, 2.6) {$\vdots$};
                \draw[thick,darkgray,rounded corners=0.5ex] (0,3) rectangle (1,3.5);
                
                \node at (2.0, 0.2) {$\cdots$};
                \node at (2.0, 1.85) {$\cdots$};
                \node at (2.0, 3.15) {$\cdots$};
                
                \begin{scope}[shift={(3.0,0)}]
                    \draw[thick,darkgray,rounded corners=0.5ex] (0,0) rectangle (1,0.5);
                    \node at (0.5, 1.15) {$\vdots$};
                    \draw[thick,darkgray,rounded corners=0.5ex] (0,1.5) rectangle (1,2);
                    \node at (0.5, 2.6) {$\vdots$};
                    \draw[thick,darkgray,rounded corners=0.5ex] (0,3) rectangle (1,3.5);
                \end{scope}
                
                \begin{scope}[shift={(4.5,0)}]
                   \draw[thick,darkgray,rounded corners=0.5ex] (0,0) rectangle (1,0.5);
                    \node at (0.5, 1.15) {$\vdots$};
                    \draw[thick,darkgray,rounded corners=0.5ex] (0,1.5) rectangle (1,2);
                \end{scope}

                \node[darkgray] at (0.5, -0.25) {\scriptsize \texttt{column $1$}};
                \node[darkgray] at (3.5, -0.25) {\scriptsize \texttt{column}};
                \node[darkgray] at (3.5, -0.5) {\scriptsize \texttt{$i-|S|$}};
                \node[darkgray] at (5, -0.25) {\scriptsize \texttt{column}};
                \node[darkgray] at (5, -0.5) {\scriptsize \texttt{$i-|S|+1$}};
                
                \node[darkgray] at (-0.6, 0.25) {\scriptsize \texttt{color $1$}};
                \node[darkgray] at (-0.6, 1.75) {\scriptsize \texttt{color $b$}};
                \node[darkgray] at (-0.6, 3.25) {\scriptsize \texttt{color $a$}};
                
                \draw[gray] (1.25, -0.25) -- (1.25, 3.65);
                \draw[gray] (2.75, -0.25) -- (2.75, 3.65);
                \draw[gray] (4.25, -0.25) -- (4.25, 3.65);
            \end{tikzpicture}
        \end{center}
        \vspace*{-0.6cm}
        \captionsetup{width=.8\linewidth}
        \caption{An arrangement of $a(i-|S|)+b$ bins into $i-|S|+1$ columns. The first $i-|S|$ columns have $a$ bins each, while the last has $b$ bins.}
        \label{figure:bins-assignment-b-neq-0}
    \end{figure}
    
    As each set $|A_S|$ is in bijection with the set of $b$-SCMs with $i+1-|S|$ segments, we conclude that
    \begin{align*}
        |Q| &= \sum_{S \subseteq \{0,1,\ldots,rp\}} (-1)^{|S|}\binom{a(i-|S|)+b}{r}^p\\
%        &= \sum_{j=0}^{rp} (-1)^j\binom{rp+1}{j}\binom{a(i-j)+b}{r}^p\\
        &= \sum_{j=0}^{i} (-1)^j\binom{rp+1}{j}\binom{a(i-j)+b}{r}^p
    \end{align*}
    as desired.
\end{proof}

Since the total number of $(p,r,a,b)$-type colored multipermutations is independently known, Theorem~\ref{theorem:main-b-neq-0} admits as a corollary a nice binomial identity, which is analytically proved in~\cite{pita-ruiz:generalized-eulerian-and-apps}.

\begin{corollary}
    Suppose $0 \leq b < a$. Then 
    \[
        \sum_{i=0}^{rp} A_{a,b,r}(p,i) = \sum_{i=0}^{rp} \sum_{j=0}^{i} (-1)^j \binom{rp+1}{j} \binom{a(i-j)+b}{r}^p = a^{rp} \frac{(rp)!}{(r!)^p}.
    \]
\end{corollary}

The definition of $A_{a,b,r}(p,i)$ given by Pita-Ruiz~\cite{pita-ruiz:generalization-eulerian,pita-ruiz:generalized-eulerian-and-apps} does not constrain $b$ to be less than $a$, as we have here. We believe that a combinatorial interpretation of the case $b \geq a$ ($b \in \mathbb{Z}$) similar to the ones we have given here could be found, although the modifications to the definition of descents --- beyond terminal and non-terminal --- seem less straight-forward. Furthermore, as Pita-Ruiz shows in \cite{pita-ruiz:generalization-eulerian}, the numbers $A_{1,4,3}(p,i)$ include negative entries, which obviously complicates any combinatorial explanation.  In this present work, we simply note that the simplification
\[
    \binom{an + b}{r}^p = \binom{a\left(n + \lfloor b/a \rfloor\right) + \left(b - a\lfloor b/a \rfloor\right)}{r}^p
\]
reduces the general case to the cases we have discussed here; when $b \geq a$,  the bins interpretation of Figure \ref{figure:bins-assignment-b-neq-0} would necessitate a final column that is taller than the others.  The algebraic manipulation above reflects the fact that, if the notion of ``weak descent'' is suitably redefined, the bins in the final column can be split into some positive number of full columns with $a$ bins each, and a final column with fewer than $a$ bins.

\section{A Combinatorial Proof of the Generalized Worpitzky Identity} %\eqref{equation:def-of-A}}
\label{sec:Worpitzky}
The ideas of the previous sections allow us to give a combinatorial derivation of the generalized Worpitzky identity in Equation~\eqref{equation:def-of-A} as well as for a summed version. (See \cite{engbers:two-comb-proofs} for the case $a=1$ and $b=0$.)

\begin{theorem}\label{theorem:sum}
Suppose $n \geq 0$, $a \geq 1$, $p \geq 1$, $r \geq 1$, and $0 \leq b < a$.  Then
\[
\sum_{k=0}^{n} \binom{ak+b}{r}^p = \sum_{i=0}^{rp} A_{a,b,r}(p,i) \binom{n+1+rp-i}{rp+1}.
\]
\end{theorem}

\begin{proof}
Consider an arrangement of bins with $n+1$ columns and $a$ bins per column, as in Figure~\ref{figure:bins-1}.

An \emph{admissible} bin placement is one formed in the following way:
\begin{enumerate}[(a)]
    \item designate one column, say column $k+1$; 
    \item amongst the $ak+b$ bins comprising the leftmost $k$ columns and the bottommost $b$ bins of column $k+1$, place $r$ copies of $1$, $r$ copies of $2$, etc.; and
    \item place one extra $1$ in bin $b+1$ of column $k+1$.  
\end{enumerate}
We remark that each admissible bin placement corresponds to an \emph{augmented $(p,r,a,b)$-type colored multipermutation}, which takes a $(p,r,a,b)$-type colored multipermutation and appends $1^{c_{b+1}}$ to the end.  We now count admissible bin placements in two ways.

The first way is quite simple. Fixing the designated column $k+1$, there are $\binom{ak+b}{r}^p$ ways to distribute $r$ copies of each of the $p$ symbols into the $ak+b$ permitted bins. Summing over all possible designated columns gives the left side of the identity.

Alternately, we may count admissible bin placements by partitioning the set of all bin placements in the following way.  We group together those bin placements whose augmented $(p,r,a,b)$-type colored multipermutations have exactly $i+1$ weak descents. Note that there is always a terminal weak descent after the final entry $1^{c_{b+1}}$ in an augmented $(p,r,a,b)$-type colored multipermutation, as the color has value larger than $b$.  Further, the augmented element $1^{c_{b+1}}$ assures that the final position of the (non-augmented) $(p,r,a,b)$-type colored multipermutation is a weak descent if and only if the color of the last entry has value strictly greater than $b$, and so by Theorem \ref{theorem:main-b-neq-0} there are $A_{a,b,r}(p,i)$ such augmented $(p,r,a,b)$-type colored multipermutations with $i+1$ weak descents.

Given a fixed augmented $(p,r,a,b)$-type colored multipermutation with exactly $i+1$ weak descents, we must enumerate the number of ways of assigning columns to the entries to recover the bin placements. Therefore, we must assign columns to each of the $rp+1$ entries of the augmented $(p,r,a,b)$-type colored multipermutation under the restriction that the columns not decrease when read from left-to-right.  Once the columns are chosen for each element, the bin is determined by the color.

We now calculate the number of ways that columns can be assigned to the elements of an augmented $(p,r,a,b)$-type colored multipermutation. Each non-terminal weak descent $i_1^{c_{j_1}} i_2^{c_{j_2}}$ requires $i_1$ and $i_2$ to lie in distinct columns, while for each ascent $i_3^{c_{j_3}} i_4^{c_{j_4}}$ the symbols $i_3$ and $i_4$ may either occur within a single column or span multiple columns. Suppose that exactly $\ell$ columns out of the total of $n+1$ contain at least one symbol. This can only be done if $\ell \geq i+1$ as each non-terminal weak descent forces a transition from one column to the next, and this further implies that exactly $\ell-(i+1)$ of the $rp+1-(i+1)=rp-i$ ascents span multiple columns. Therefore, the number of such column assignments is
\[
    \sum_{\ell=i+1}^{n+1} \binom{n+1}{\ell}\binom{rp-i}{\ell-i-1} = \sum_{\ell=i}^{n+1} \binom{n+1}{\ell}\binom{rp-i}{rp-\ell+1} = \binom{n+1+rp-i}{rp+1}.
\]
Each choice of columns corresponds to a unique admissible bin placement, and each admissible bin placement is represented by such a choice of columns. This gives the right side of the identity.
\end{proof}

A combinatorial proof of the generalization of Worpitzky's identity from \eqref{equation:def-of-A} follows quickly from the ideas contained in the proof of Theorem~\ref{theorem:sum}. 

\begin{porism}
Suppose $n \geq 0$, $a \geq 1$, $p \geq 1$, $r \geq 1$, and $0 \leq b < a$.  Then
\[
\binom{an+b}{r}^{p} = \sum_{i=0}^{rp} A_{a,b,r}(p,i) \binom{n+rp-i}{rp}.
\]
\end{porism}

\begin{proof}
Consider only the admissible bin placements from the proof of Theorem~\ref{theorem:sum} in which the designated column is $n+1$. The number of ways to assign $r$ copies of each of the $p$ symbols into the $an+b$ total boxes is $\binom{an+b}{r}^p$, the left side of the identity.

As before, we find the right side of the identity by counting the number of ways that columns can be assigned to each of the entries. Here we must add the condition that the final column must be one of those that is selected. Thus, in this case, the number of such column assignments is
\[
  \sum_{\ell=i+1}^n \binom{n}{\ell-1}\binom{rp-i}{\ell-i-1} = \sum_{\ell=i+1}^n \binom{n}{\ell-1}\binom{rp-i}{rp-\ell+1} = \binom{n + rp - i}{rp}.
\]
The remaining details are identical to the proof of Theorem~\ref{theorem:sum}.
\end{proof}

At this point, it is worth mentioning that the form of the generalization of Worpitzky's identity in Equation~\eqref{equation:def-of-A} is motivated by $f$-Eulerian numbers \cite{pita-ruiz:generalization-eulerian,stanleyec1}. An alternate way to generalize Worpitzky's identity, fitting work from \cite{engbers:two-comb-proofs}, would be to define numbers $A'_{a,b,r}(p,i)$ by
\[
    \binom{an+b}{r}^p = \sum_{i=0}^{rp} A'_{a,b,r}(p,i) \binom{n+i}{rp}.
\]
The translation here is $A'_{a,b,r}(p,i) = A_{a,b,r}(p,rp-i)$, and so this would lead to an interpretation of $A'_{a,b,r}(p,i)$  involving ``ascents'' instead of ``weak descents.''  We leave the details to the interested reader.

%=================================================================================
%References
%=================================================================================

%\bibliography{bib}
\bibliography{bib.bib}
\bibliographystyle{plain}

\end{document}